\documentclass{amsart}
\usepackage{amsmath,amssymb,amsthm,verbatim}
\newtheorem{theorem}{Theorem}
\newtheorem{lemma}{Lemma}

\newtheorem{proposition}{Proposition}
\theoremstyle{definition}

\begin{document}
\title{An application of sums of triple products of binomials}
\author{Michael~J.~J.~Barry}
\curraddr{15 River Street Apt 205\\
Boston,  MA 02108}
\email{mbarry@allegheny.edu}
\begin{abstract}
We prove that a certain family of sums of products of three binomials has alternating behavior modulo a prime $p$.  To accomplish this we rewrite these sums as signed sums of  products of three binomials, the better to handle $p$,  and we give closed-form expressions for two related sums of signed products of three binomials.
\end{abstract}
\maketitle

\section{Statement of Main Result}
\begin{theorem}\label{Main}
Suppose that $p$ is a prime,  and that $k$, $c$, and $d$ are integers satisfying $ 1 \leq k \leq c \leq d <c+d \leq p$.   Define the functions $C_{p,c,d,k}$ and $D_{p,c,d,k}$ for every integer $\ell \in [1,c+d+1-k]$ by
\[C_{p,c,d,k}(\ell)=\sum_{j=1}^{c+1-k} \binom{k+j-2}{k-1 }\binom{c+d-k }{d+j-1 }\binom{p-c-d+2k-2}{k+j-1-\ell}\]
and
\[D_{p,c,d,k}(\ell)=\sum_{j=1}^{d+1-k}\binom{d-j}{k-1}\binom{ c+d-k}{ j-1}\binom{p-c-d+2k-2}{p+k+j-d-1-\ell},\]
and let $f_{p,c,d,k}=C_{p,c,d,k}+(-1)^k D_{p,c,d,k}$.  Then $f_{p,c,d,k}(1) \not \equiv 0 \mod p$ and $f_{p,c,d,k}(\ell) \equiv (-1)^{\ell-1} f_{p,c,d,k}(1) \mod p$ for every integer $\ell \in [1,c+d+1-k]$.
\end{theorem}

To save space we will write $f$ for $f_{p,c,d,k}$, $C$ for $C_{p,c,d,k}$,  and $D$ for $D_{p,c,d,k}$. 

Our first step in proving Theorem~\ref{Main}
will be to rewrite $C(\ell)$ and $D(\ell)$ as
\[
\sum_{r=0}^{k-1}(-1)^r\binom{c-1-r}{k-1-r}\binom{c+d-k}{r}\binom{p+k-r-2}{c-\ell-r}
\]
and
\[
\sum_{j=0}^{k-1}(-1)^j\binom{d-1-j}{k-1-j}\binom{c+d-k}{j}\binom{p+k-2-j}{\ell+k-2-c-j},\]
respectively.  This will allow us to deal with the prime $p$ more effectively.

In dealing with $p$, we will give in Lemma~\ref{FandG} closed-form evaluations of two sums of signed products of three binomials:
\[\sum_{r=0}^{k-1}(-1)^r\binom{c-1-r}{k-1-r}\binom{c+d-k}{r}\binom{k-r-1}{k-1-j}\]
and
\[\sum_{r=0}^{k-1}(-1)^r\binom{d-1-r}{k-1-r}\binom{c+d-k}{r}\binom{k-1-r}{j}.
\]

In Section~\ref{Motivation}, we indicate how this result arose; in Sections~\ref{CEval} and~\ref{DEval}, we derive the rewrites of $C(\ell)$ and $D(\ell)$; in Section~\ref{MainResult}, we proof Theorem~\ref{Main}; and in Section~\ref{Others}, we suggest possible related results.
\section{Motivation}\label{Motivation}

How does the function $f=f_{p,c,d,k}$ arise?

Let $p$ be a prime number, $F$ a field of characteristic $p$, and $G$ a cyclic group of order $q=p^a>1$.  Up to isomorphism, there is a unique indecomposable $F G$-module $V_q$ of dimension $q$~\cite[pp. 24--25]{Alperin}.  Let $g$ be a generator of $G$.  Then there is an ordered $F$-basis $(v_1,v_2,\dots,v_q)$ of $V_q$ such that $g v_1=v_1$ and $g v_i=v_{i-1}+v_i$ if $i>1$, that is, the matrix of $g$ with respect to this basis is a full Jordan block of eigenvalue $1$.  For an integer $i \in [1,q]$ define the vector space $V_i$ over $F$ by $V_i=\langle v_1,\dots,v_i \rangle$.  Then $V_i$ is an indecomposable $F G$-module and $\{V_1,\dots,V_q\}$ is a complete set of indecomposable $F G$-modules~\cite[pp. 24--25]{Alperin}.  For integers $m$ and $n$ in $[1,q]$, $\{v_i \otimes v_j \mid 1 \leq i \leq m, 1 \leq j \leq n\}$ is an $F$-basis of $V_m \otimes V_n$.  But $\mathcal{B}=\{v_{i,j}=v_i \otimes g^{n-i}v_j \mid 1 \leq i \leq m, 1 \leq j \leq n\}$ is another basis with the nice property that $(g -1)(v_{i,j})=v_{i-1,j}+v_{i,j-1}$~\cite[Lemma 1]{B2017}.

Specialize $m$ and $n$ to $m=p+c$ and $n=p+d$ where $1 \leq c \leq d <c+d \leq p$.  Then by~\cite[Theorem 1]{Renaud1979}, 
\[V_c \otimes V_d \cong \bigoplus_{k=1}^c V_{\lambda_k}\]
where $\lambda_k=c+d-2k+1$.  And by repeated applications of~\cite[Corollary 1]{B2011_0} or by~\cite[Theorem 2]{Renaud1979} 
\begin{align*}
V_{p+c}\otimes V_{p+d}&\cong \bigoplus_{k=1}^cV_{2p+\lambda_k}
\oplus (d-c)\dot V_{2p} \oplus \bigoplus_{k=1}^c V_{2p-\lambda_k}\oplus V_{p-d}\otimes V_{p-c}\\
&\cong \bigoplus_{k=1}^cV_{2p+\lambda_k}
\oplus (d-c)\dot V_{2p} \oplus \bigoplus_{k=1}^c V_{2p-\lambda_k} \oplus (p-c-d) \cdot V_p
\oplus V_c \otimes V_d\\
&\cong \bigoplus_{k=1}^cV_{2p+\lambda_k}
\oplus (d-c)\dot V_{2p} \oplus \bigoplus_{k=1}^c V_{2p-\lambda_k} \oplus (p-c-d) \cdot V_p
\oplus \bigoplus_{k=1}^c V_{\lambda_k}.
\end{align*}

The function $f$ arises in identifying a generator for the cyclic module $V_{2p-\lambda_k}$, where $1 \leq k \leq c$,  in terms of the basis $\mathcal{B}$ as we now explain.  Define $y_{c+d+1-k}\in V_{p+c} \otimes V_{p+d}$ by
\begin{align*}
y_{c+d+1-k}&=\sum_{j=1}^{c+1-k} \binom{k+j-2}{k-1 }\binom{c+d-k }{d+j-1 }v_{p+k+j-1,p+1-j}\\
&\quad+(-1)^k \sum_{j=1}^{d+1-k}\binom{d-j}{k-1}\binom{ c+d-k}{ j-1}v_{p+k+j-d-1,p+d+1-j}.
\end{align*}

The coefficients in $y_{c+d+1-k}$ come from a $p \times p$ matrix $B(m,n;p)$ with $(i,j)$-entry 
\[\binom{c-i}{c+d-i-j+1}\binom{i+j-2}{i-1} \in F\] defined in Norman~\cite[p. 431]{Norman2008}.

The $(c+d+2-k-j,j)$-entry of $B(c,d;p)$, where $1 \leq j \leq d+1-k$, is
\[\binom{k+j-d-2}{k-1}\binom{c+d-k}{c+d+1-k-j}=(-1)^{k-1}\binom{d-j}{k-1}\binom{c+d-k}{j-1}\]
and the $(c+2-k-j,d+j)$-entry of $B(c,d;p)$, where $1 \leq j \leq c+1-k$, is
\[\binom{k+j-2}{k-1}\binom{c+d-k}{c+1-k-j}=\binom{k+j-2}{k-1}\binom{c+d-k}{d+j-1}.
\]

Thus the coefficients are from the $(c+d+1-k)$-th anti-diagonal of $B(c,d;p)$,  except that the first $d+1-k$ coefficients have opposite sign.

The coefficient of $v_{\ell,c+d+2-k-\ell}$ in $(g-1)^{2p-\lambda_{k}-1}(y_{c+d+1-k})$ is
\begin{align*}
&\sum_{j=1}^{c+1-k} \binom{k+j-2}{k-1 }\binom{c+d-k }{d+j-1 }\binom{2p-\lambda_k-1}{p+k+j-\ell-1}\\
&\quad+(-1)^k \sum_{j=1}^{d+1-k}\binom{d-j}{k-1}\binom{ c+d-k}{ j-1}\binom{2p-\lambda_k-1}{p+k+j-d-1-\ell}\\
&=\sum_{j=1}^{c+1-k} \binom{k+j-2}{k-1 }\binom{c+d-k }{d+j-1 }\binom{2p-c-d+2k-2}{p+k+j-\ell-1}\\
&\quad+(-1)^k \sum_{j=1}^{d+1-k}\binom{d-j}{k-1}\binom{ c+d-k}{ j-1}\binom{2p-c-d+2k-2}{p+k+j-d-1-\ell}\\
&= \sum_{j=1}^{c+1-k} \binom{k+j-2}{k-1 }\binom{c+d-k }{d+j-1 }\binom{p-c-d+2k-2}{k+j-\ell-1}\\
&\quad+(-1)^k \sum_{j=1}^{d+1-k}\binom{d-j}{k-1}\binom{ c+d-k}{ j-1}\binom{p-c-d+2k-2}{k+j-d-1-\ell}\\
\end{align*}
since $\text{char }F=p$.  But this equals $(-1)^{\ell-1}f(1)$ in $F$ by Theorem ~\ref{Main}.  Thus
\[(g-1)^{2p-\lambda_{k}-1}(y_{c+d+1-k})=f(1)\sum_{\ell=1}^{c+d+1-k}(-1)^{\ell-1}v_{\ell,c+d+1-k-\ell}.\]
This equation and the following easily verified equation
\[(g-1)\left(\sum_{\ell=1}^{c+d+1-k}(-1)^{\ell-1}v_{\ell,c+d+1-k-\ell}\right)=0\]
 show that $y_{c+d+1-k}$ generates a cyclic indecomposable module of dimension $2p-\lambda_{k}$.
 
Though the material in this section and that of Norman~\cite{Norman2008} are clearly related,  our use of $B(c,d;p)$ and Norman's use, for example, in~\cite[Lemma 11]{Norman2008} are different.

\section{Evaluating $C$}\label{CEval}
\begin{proposition}\label{CProposition}
For $\ell \in [1,c+d+1-k]$,
\[C(\ell)=\sum_{r=0}^{k-1}(-1)^r\binom{c-1-r}{k-1-r}\binom{c+d-k}{r}\binom{p+k-r-2}{c-\ell-r}.\]
\end{proposition}
The proof of Proposition~\ref{CProposition} will require the following two results.
\begin{lemma}\label{CLemma1}
For every integer $j \in[0,c-k]$,
\[\binom{c-1}{k-1}-\binom{k-1+j}{k-1}=\sum_{r=1}^{k-1}(-1)^{r-1}\binom{c-k-j}{r}\binom{c-1-r}{k-1-r}.\]
\end{lemma}
\begin{proof}
Note that this is equivalent to proving
\[\binom{k-1+j}{k-1}=\sum_{r=0}^{k-1}(-1)^{r}\binom{c-k-j}{r}\binom{c-1-r}{k-1-r}.\]
Now
\begin{align*}
\sum_{r=0}^{k-1}&(-1)^{r}\binom{c-k-j}{r}\binom{c-1-r}{k-1-r}\\
&=\sum_{r=0}^{k-1}(-1)^{r}\binom{c-k-j}{r}\binom{c-1-r}{c-k}\\
&=(-1)^{c-1+c-k}\binom{c-k-j-(c-k)-1}{c-1-(c-k)-0}&& \text{by~\cite[Equation (5.25)]{GKP}}\\
&=(-1)^{k+1}\binom{-j-1}{k-1}\\
&=(-1)^{k+1}(-1)^{k-1}\binom{k-1+j+1-1}{k-1}\\
&=\binom{k-1+j}{k-1}.
\end{align*}
\end{proof}

\begin{lemma}\label{L-KToMinus1}
\[\sum_{r=0}^{k-1}(-1)^{r+1}\binom{c-1-r}{k-1-r}\binom{c+d-k}{r}\left(\sum_{j=\ell-k}^{-1}\binom{c+d-k-r}{c-k-j-r}\binom{p-c-d+2k-2}{k+j-\ell}\right)=0\]
for every integer $\ell \in [1-k,-1]$.
\end{lemma}
\begin{proof}
It suffices to show that
\[\sum_{r=0}^{k-1}(-1)^{r+1}\binom{c-1-r}{k-1-r}\binom{c+d-k}{r}\binom{c+d-k-r}{c-k-j-r}=0
\]
for every $j \in [1-k,-1]$.

Now
\begin{align*}
\sum_{r=0}^{k-1}&(-1)^{r+1}\binom{c-1-r}{k-1-r}\binom{c+d-k}{r}\binom{c+d-k-r}{c-k-j-r}\\
&=\sum_{r=0}^{k-1}(-1)^{r+1}\binom{c-1-r}{k-1-r}\binom{c+d-k}{c-k-j}\binom{c-k-j}{r}\\
&=-\binom{c+d-k}{c-k-j}\left(\sum_{r=0}^{k-1}(-1)^{r}\binom{c-1-r}{k-1-r}\binom{c-k-j}{r}\right)\\
&=-\binom{c+d-k}{c-k-j}\binom{-j-1}{k-1}, \qquad \text{by~\cite[Equation (5.25)]{GKP}.}\\
\end{align*}
But $\binom{-j-1}{k-1}=0$ for $j=1-k, 2-k,\dots,-1$.
\end{proof}

\begin{proof}[Proof of Proposition~\ref{CProposition}]
First note that for an integer $\ell \in [1,c+d+1-k]$,
\begin{align*}
C(\ell)&=\sum_{j=0}^{c-k}\binom{k+j-1}{k-1}\binom{c+d-k}{d+j}\binom{p-c-d+2k-2}{k+j-\ell}\\
&=\binom{k-1}{k-1}\binom{c+d-k}{d}\binom{p-c-d+2k-2}{k-\ell}\\&+\quad
\sum_{j=1}^{c-k}\binom{k+j-1}{k-1}\binom{c+d-k}{d+j}\binom{p-c-d+2k-2}{k+j-\ell}.
\end{align*}
Letting $g(\alpha,\beta)=\sum_{r=\alpha}^\beta \binom{c+d-k}{d+r}\binom{p-c-d+2k-2}{k+r-\ell}$, we see that
\begin{align*}
C(\ell)&=\binom{k-1}{k-1}g(0,c-k)+\sum_{j=1}^{c-k}\left(\binom{k-1+j}{k-1}-\binom{k-1}{k-1}\right)\binom{c+d-k}{d+j}\binom{p-c-d+2k-2}{k+j-\ell}.
\end{align*}
Since $\binom{k}{k-1}-\binom{k-1}{k-1}=\binom{k-1}{k-2}$,
\begin{align*}
C(\ell)&=\binom{k-1}{k-1}g(0,c-k)+\binom{k-1}{k-2}\sum_{j=1}^{c-k}\binom{c+d-k}{d+j}\binom{p-c-d+2k-2}{k+j-\ell}\\
&\quad+\sum_{j=2}^{c-k}\left(\binom{k+j-1}{k-1}-\binom{k-1}{k-2}-\binom{k-1}{k-1}\right)\binom{c+d-k}{d+j}\binom{p-c-d+2k-2}{k+j-\ell}\\
&=\binom{k-1}{k-1}g(0,c-k)+\binom{k-1}{k-2}\sum_{j=1}^{c-k}\binom{c+d-k}{d+j}\binom{p-c-d+2k-2}{k+j-\ell}\\
&\quad+\sum_{j=2}^{c-k}\left(\binom{k+j-1}{k-1}-\binom{k}{k-1}\right)\binom{c+d-k}{d+j}\binom{p-c-d+2k-2}{k+j-\ell}\\
&=\binom{k-1}{k-1}g(0,c-k)+\binom{k-1}{k-2}g(1,c-k)\\
&\quad+\sum_{j=2}^{c-k}\left(\binom{k+j-1}{k-1}-\binom{k}{k-1}\right)\binom{c+d-k}{d+j}\binom{p-c-d+2k-2}{k+j-\ell}.\\
\end{align*}
Continuing is this way, we get
\[
C(\ell)=\binom{k-1}{k-1}g(0,c-k)+\sum_{j=1}^{c-k}\binom{k-2+j}{k-2}g(j,c-k).
\]
By~\cite[Equation (5.23)]{GKP}
\[\sum_{r} \binom{c+d-k}{d+r}\binom{p-c-d+2k-2}{k+r-\ell}=\binom{p+k-2}{c+d-k-d+k-\ell}=\binom{p+k-2}{c-\ell}.\]
But
\[\sum_{r} \binom{c+d-k}{d+r}\binom{p-c-d+2k-2}{k+r-\ell}=\sum_{r=\ell-k}^{c-k} \binom{c+d-k}{d+r}\binom{p-c-d+2k-2}{k+r-\ell}.\]
Thus
\[g(0,c-k)=\binom{p+k-2}{c+d-k-d+k-\ell}-g(\ell-k,-1),\]
and
\begin{align*}
g(j,c-k)&=\sum_{r=j}^{c-k} \binom{c+d-k}{d+r}\binom{p-c-d+2k-2}{k+r-\ell}\\
&=g(0,,c-k)-g(\ell-k,j-1)\\
&=\binom{p+k-2}{c-\ell}-g(\ell-k,-1)-g(0,j-1).\\.
\end{align*}

Hence
\begin{align*}
C(\ell)&=\binom{k-1}{k-1}\left(\binom{p+k-2}{c-\ell}-g(\ell-k,-1)\right)\\
&\quad+\sum_{j=1}^{c-k}\binom{k-2+j}{k-2}\left(\binom{p+k-2}{c-\ell}-g(\ell-k,-1)-g(0,j-1)\right).
\end{align*}
But
\[\binom{k-1}{k-1}+\sum_{j=1}^{c-k}\binom{k-2+j}{k-2}=\binom{c-1}{k-1}.
\]
Thus
\begin{align*}
C(\ell)&=\binom{c-1}{k-1}\left(\binom{p+k-2}{c-\ell}-g(\ell-k,-1)\right)
-\sum_{j=1}^{c-k}\binom{k-2+j}{k-2}g(0,j-1).
\end{align*}
Denote $\sum_{j=1}^{c-k}\binom{k-2+j}{k-2}g(0,j-1)$ by $H(\ell)$.  Then
\begin{align*}
H(\ell)&=\sum_{j=1}^{c-k}\binom{k-2+j}{k-2}\left(\sum_{r=0}^{j-1}\binom{c+d-k}{d+r}\binom{p-c-d+2k-2}{k+r-\ell}\right)\\
&=\sum_{r=0}^{c-k-1}\binom{c+d-k}{c-k-r}\binom{p-c-d+2k-2}{k+r-\ell}
\left(\sum_{j=r}^{c-k-1}\binom{k-1+j}{k-2}  \right)\\
&=\sum_{r=0}^{c-k-1}\binom{c+d-k}{c-k-r}\binom{p-c-d+2k-2}{k+r-\ell}
\left(\binom{c-1}{k-1}-\binom{k-1+r}{k-1}\right)\\
&=\sum_{r=0}^{c-k-1}\binom{c+d-k}{c-k-r}\binom{p-c-d+2k-2}{k+r-\ell}
\left(\sum_{z=1}^{k-1}(-1)^z \binom{c-k-r}{z}\binom{c-1-z}{k-1-z}\right)
\end{align*}
by Lemma~\ref{CLemma1}.  Continuing
\begin{align*}
H(\ell)&=\sum_{z=1}^{k-1}(-1)^z\binom{c-1-z}{k-1-z}
\left(\sum_{r=0}^{c-k-1}\binom{c+d-k}{c-k-r}\binom{c-k-r}{z}\binom{p-c-d+2k-2}{k+r-\ell}\right)\\
&=\sum_{z=1}^{k-1}(-1)^z\binom{c-1-z}{k-1-z}
\left(\sum_{r=0}^{c-k-1}\binom{c+d-k}{z}\binom{c+d-k-z}{c-k-r-z}\binom{p-c-d+2k-2}{k+r-\ell}\right)\\
&=\sum_{z=1}^{k-1}(-1)^z\binom{c-1-z}{k-1-z}\binom{c+d-k}{z}
\left(\sum_{r=0}^{c-k-1}\binom{c+d-k-z}{c-k-r-z}\binom{p-c-d+2k-2}{k+r-\ell}\right)\\
&=\sum_{z=1}^{k-1}(-1)^z\binom{c-1-z}{k-1-z}\binom{c+d-k}{z}
\left(\binom{p+k-z-2}{c-\ell-z}-\sum_{r=\ell-k}^{-1}\binom{c+d-k-z}{c-k-r-z}\binom{p-c-d+2k-2}{k+r-\ell}\right)\\
\end{align*}
by~\cite[Equation (5.23)]{GKP}.

Thus
\[C(\ell)=\sum_{z=0}^{k-1}(-1)^z\binom{c-1-z}{k-1-z}\binom{c+d-k}{z}\binom{p+k-z-2}{c-\ell-z}\]
since by Lemma~\ref{L-KToMinus1},
\[\sum_{z=0}^{k-1}(-1)^z\binom{c-1-z}{k-1-z}\binom{c+d-k}{z}\left(
\sum_{r=\ell-k}^{-1} \binom{c+d-k-z}{c-k-r-z}\binom{p-c-d+2k-2}{k+r-\ell}\right)=0.\]
\end{proof}

\section{Evaluating $D$}\label{DEval}
\begin{proposition}\label{DProposition}
For $\ell \in [1,c+d+1-k]$,
\[D(\ell)=\sum_{j=0}^{k-1}(-1)^j\binom{d-1-j}{k-1-j}\binom{c+d-k}{j}\binom{p+k-2-j}{\ell+k-2-c-j}.\]
\end{proposition}
The proof of Proposition~\ref{DProposition} will require the following result whose proof is similar to the proof of Lemma~\ref{CLemma1}.
\begin{lemma}\label{DLemma1}
When $ 0 \leq j \leq k-1$,
\[\binom{d-1}{k-1}-\binom{d-1-j}{k-1}=\sum_{r=1}^{k-1}(-1)^{r-1}\binom{d-1-r}{k-1-r}\binom{j}{r}.\]
\end{lemma}

\begin{proof}[Proof of Proposition~\ref{DProposition}]
First note that for an integer $\ell \in [1,c+d+1-k]$,
\[D(\ell)=\sum_{j=0}^{\ell+k-2-c}\binom{d-1-j}{k-1}\binom{ c+d-k}{ j}\binom{p-c-d+2k-2}{\ell+k-2-c-j}.
\]
Then
\begin{align*}
D(\ell)&=\sum_{j=0}^{\ell+k-3-c}\binom{d-1-j}{k-1}\binom{ c+d-k}{ j}\binom{p-c-d+2k-2}{\ell+k-2-c-j}\\
&\quad+\binom{c+d+1-k-\ell}{k-1}\binom{c+d-k}{\ell+k-2-c}\binom{p-c-d+2k-2}{0}.
\end{align*}
Letting $h(\alpha,\beta)=\sum_{j=\alpha}^{\beta}\binom{ c+d-k}{ j}\binom{p-c-d+2k-2}{\ell+k-2-c-j}$ and noting that $\binom{d-1-(\ell+k-2-c)}{k-1}=\binom{c+d+1-k-\ell}{k-1}$,
\begin{align*}
D(\ell)&=\binom{c+d+1-k-\ell}{k-1}h(0,\ell+k-2-c)\\
&\quad+\sum_{r=0}^{\ell+k-c-3}\left(\binom{d-1-r}{k-1}- \binom{c+d+1-k-\ell}{k-1}\right)h(r,\ell+k-3-c)\\
&=\binom{c+d+1-k-\ell}{k-1}h(0,\ell+k-2-c)\\
&\quad+\sum_{r=0}^{\ell+k-c-4}\left(\binom{d-1-r}{k-1}- \binom{c+d+1-k-\ell}{k-1}\right)h(r,\ell+k-3-c)\\
&\quad+\binom{c+d+1-k-\ell}{k-2}\binom{c+d-j}{\ell+k-c-3}\binom{p-c-d+2k-2}{1}\\
&=\binom{c+d+1-k-\ell}{k-1}h(0,\ell+k-2-c)\\
&\quad+\binom{c+d+1-k-\ell}{k-2}h(0,\ell+k-3-c)\\
&\quad+\sum_{r=0}^{\ell+k-c-4}\left(\binom{d-1-r}{k-1}- \binom{c+d+2-k-\ell}{k-1}\right)
h(r,\ell+k-3-c)\\
\end{align*}
using the fact that $\binom{c+d+1-k-\ell}{k-1}+\binom{c+d+1-k-\ell}{k-2}=\binom{c+d+2-k-\ell}{k-1}$.

Continuing in this way we get
\begin{align*}
D(\ell)&=\binom{c+d+1-k-\ell}{k-1}h(0,\ell+k-2-c)\\
&\quad+\sum_{r=0}^{\ell+k-3-c}\binom{c+d+1-k-\ell+r}{k-2}h(0,\ell+k-3-c-r)\\
&=\binom{c+d+1-k-\ell}{k-1}h(0,\ell+k-2-c)\\
&\quad+\sum_{r=0}^{\ell+k-3-c}\binom{c+d+1-k-\ell+r}{k-2}\left(h(0,\ell+k-2-c)-h(\ell+k-2-c-r,\ell+k-2-c)\right).
\end{align*}

Hence
\begin{align*}
D(\ell)&=h(0,\ell+k-2-c)\left(\binom{c+d+1-k-\ell}{k-1}+\sum_{r=0}^{\ell+k-3-c}\binom{c+d+1-k-\ell+r}{k-2}\right)\\
&\quad-\sum_{r=0}^{\ell+k-3-c}\binom{c+d+1-k-\ell+r}{k-2}h(\ell+k-2-c-r,\ell+k-2-c).
\end{align*}
But
\[h(0,\ell+k-2-c)=\sum_{j=0}^{\ell+k-2-c}\binom{ c+d-k}{ j}\binom{p-c-d+2k-2}{\ell+k-2-c-j}=\binom{p+k-2}{\ell+k-c-2}\]
by Vandermonde's convolution~\cite[Equation (5.22)]{GKP}  and
\[\binom{c+d+1-k-\ell}{k-1}+\sum_{r=0}^{\ell+k-3-c}\binom{c+d+1-k-\ell+r}{k-2}=\binom{d-1}{k-1}.
\]

Thus
\begin{align*}
D(\ell)&=\binom{d-1}{k-1}\binom{p+k-2}{\ell+k-c-2}\\
&\quad-\sum_{r=0}^{\ell+k-3-c}\binom{c+d+1-k-\ell+r}{k-2}h(\ell+k-2-c-r,\ell+k-2-c).
\end{align*}
Denote $\sum_{r=0}^{\ell+k-3-c}\binom{c+d+1-k-\ell+r}{k-2}h(\ell+k-2-c-r,\ell+k-2-c)$ by $E(\ell)$.
Then
\begin{align*}
E(\ell)&=\sum_{r=0}^{\ell+k-3-c}\binom{c+d+1-k-\ell+r}{k-2}\left( \sum_{j=\ell+k-2-c-r}^{\ell+k-2-c}\binom{ c+d-k}{ j}\binom{p-c-d+2k-2}{\ell+k-2-c-j}\right)\\
&=\sum_{j=1}^{\ell+k-2-c}\binom{ c+d-k}{ j}\binom{p-c-d+2k-2}{\ell+k-2-c-j}
\left( \sum_{r=\ell+k-2-c-j}^{\ell+k-3-c}\binom{c+d+1-k-\ell+r}{k-2}\right)\\
&=\sum_{j=1}^{\ell+k-2-c}\binom{ c+d-k}{ j}\binom{p-c-d+2k-2}{\ell+k-2-c-j}
\left(\sum_{r=0}^{j-1}\binom{d-2-r}{k-2}\right)\\
&=\sum_{j=1}^{\ell+k-2-c}\binom{ c+d-k}{ j}\binom{p-c-d+2k-2}{\ell+k-2-c-j}
\left(\binom{d-1}{k-1}-\binom{d-j-1}{k-1}\right)\\
&=\sum_{j=1}^{\ell+k-2-c}\binom{ c+d-k}{ j}\binom{p-c-d+2k-2}{\ell+k-2-c-j}
\left(\sum_{r=1}^{k-1}(-1)^{r-1}\binom{d-1-r}{k-1-r}\binom{j}{r}\right)
\end{align*}
by Lemma~\ref{DLemma1}.
So
\begin{align*}
E(\ell)&=\sum_{r=1}^{k-1}(-1)^{r-1}\binom{d-1-r}{k-1-r}\left(\sum_{j=1}^{\ell+k-2-c}    \binom{j}{r}\binom{ c+d-k}{ j}\binom{p-c-d+2k-2}{\ell+k-2-c-j}\right)\\
&=\sum_{r=1}^{k-1}(-1)^{r-1}\binom{d-1-r}{k-1-r}\binom{c+d-k}{r}
\left(\sum_{j=1}^{\ell+k-2-c}\binom{c+d-k-r}{j-r}\binom{p-c-d+2k-2}{\ell+k-2-c-j}\right)\\
&=\sum_{r=1}^{k-1}(-1)^{r-1}\binom{d-1-r}{k-1-r}\binom{c+d-k}{r}
\left(\sum_{j=r}^{\ell+k-2-c}\binom{c+d-k-r}{j-r}\binom{p-c-d+2k-2}{\ell+k-2-c-j}\right)\\
&=\sum_{r=1}^{k-1}(-1)^{r-1}\binom{d-1-r}{k-1-r}\binom{c+d-k}{r}
\left(\sum_{j=0}^{\ell+k-2-c-r}\binom{c+d-k-r}{j}\binom{p-c-d+2k-2}{\ell+k-2-c-j-r}\right)\\
&=\sum_{r=1}^{k-1}(-1)^{r-1}\binom{d-1-r}{k-1-r}\binom{c+d-k}{r}\binom{p+k-2-r}{\ell+k-c-2-r}
\end{align*}
by Vandermonde's convolution~\cite[Equation (5.22)]{GKP}.
Thus
\begin{align*}
D(\ell)&=\binom{d-1}{k-1}\binom{p+k-2}{\ell+k-c-2}-E(\ell)\\
&=\sum_{r=0}^{k-1}(-1)^{r}\binom{d-1-r}{k-1-r}\binom{c+d-k}{r}\binom{p+k-2-r}{\ell+k-c-2-r}.
\end{align*}
\end{proof}

\section{Proof of Main Result}\label{MainResult}
In our proof of the Theorem~\ref{Main}, we will use the following result.
\begin{lemma}\label{FandG}
Let $c$, $d$, and $k$ be positive integers such that $k \leq c \leq d$.  For every integer $\ell \in[1,c+d-k]$, define 
\[F_{c,d,k}(\ell)=\sum_{r=0}^{k-1}(-1)^r\binom{c-1-r}{k-1-r}\binom{c+d-k}{r}\binom{k-1-r}{k+\ell-c-1}\]
and
\[G_{c,d,k}(\ell)=\sum_{r=0}^{k-1}(-1)^r\binom{d-1-r}{k-1-r}\binom{c+d-k}{r}\binom{k-1-r}{c-\ell}.\]
Then $F_{c,d,k}(\ell)=(-1)^{\ell-c}\binom{\ell-1}{c-k}\binom{c+d-k-\ell}{d-k}$ and $G_{c,d,k}(\ell)=(-1)^{\ell+k-c-1}\binom{\ell-1}{c-k}\binom{c+d-k-\ell}{d-k}
$.  Hence $F_{c,d,k}(\ell)+(-1)^k G_{c,d,k}(\ell)=0$ for every integer $\ell \in[1,c+d-k]$.
\end{lemma}
\begin{proof}
We will only prove the result for $F_{c,d,k}$ as the proof for $G_{c,d,k}$ is similar.
We must show that for every integer $\ell \in [1,c+d-k]$,
\[\sum_{r=0}^{k-1}(-1)^r\binom{c-1-r}{k-1-r}\binom{c+d-k}{r}\binom{k-1-r}{k+\ell-c-1}
=(-1)^{\ell-c}\binom{\ell-1}{c-k}\binom{c+d-k-\ell}{d-k}.\]
Replacing $c-\ell$ by $j$, it suffices to  show
\[\sum_{r=0}^{k-1}(-1)^r\binom{c-1-r}{k-1-r}\binom{c+d-k}{r}\binom{k-1-r}{k-1-j}=
(-1)^j\binom{c-j-1}{c-k}\binom{d-k+j}{d-k}.
\]
Now
\begin{align*}
\sum_{r=0}^{k-1}&(-1)^r\binom{c-1-r}{k-1-r}\binom{c+d-k}{r}\binom{k-1-r}{k-1-j}\\
&=\sum_{r=0}^{j}(-1)^r\binom{c-1-r}{k-1-r}\binom{c+d-k}{r}\binom{k-1-r}{k-1-j}.
\end{align*}
But
\[\binom{c-1-r}{k-1-r}\binom{k-1-r}{k-1-j}=\binom{c-j-1}{c-k}\binom{c-1-r}{j-r}.\]
Hence
\begin{align*}
\sum_{r=0}^{j}&(-1)^r\binom{c-1-r}{k-1-r}\binom{c+d-k}{r}\binom{k-1-r}{k-1-j}\\
&=\sum_{r=0}^{j}(-1)^r\binom{c-j-1}{c-k}\binom{c-1-r}{j-r}\binom{c+d-k}{r}\\
&=\binom{c-j-1}{c-k}\sum_{r=0}^{j}(-1)^r\binom{c-1-r}{c-j-1}\binom{c+d-k}{r}\\
&=\binom{c-j-1}{c-k}(-1)^{c-1+(c-j-1)}\binom{c+d-k-(c-j-1)-1}{c-1-(c-j-1)} && \text{by~\cite[Equation (5.25)]{GKP}}\\
&=(-1)^j \binom{c-j-1}{c-k}\binom{d-k+j}{j}\\
&=(-1)^j \binom{c-j-1}{c-k}\binom{d-k+j}{d-k}.
\end{align*}
\end{proof}
\begin{proof}[Proof of Theorem~\ref{Main}]
First  if $r<k-1$, then
\[\binom{p+k-r-2}{c-1-r}=\frac{(p+k-r-2)!}{(c-1-r)!(p+k-c-1)!}\equiv 0 \mod p\]
since $k \leq c$.  Thus
\begin{align*}
C(1)&=\sum_{r=0}^{k-1}(-1)^r\binom{c-1-r}{k-1-r}\binom{c+d-k}{r}\binom{p+k-r-2}{c-1-r}\\
& \equiv(-1)^{k-1}\binom{c-k}{0}\binom{c+d-k}{k-1}\binom{p-1}{c-k} \mod p\\
& \not \equiv 0 \mod p.
\end{align*}

Second
\[
D(1)=\sum_{j=0}^{k-1}(-1)^j\binom{d-1-j}{k-1-j}\binom{c+d-k}{j}\binom{p+k-2-j}{k-1-c-j}
=0
\]
since $k-1-c-j<0$.  Hence $f(1) \not \equiv 0 \mod p$.

For $\ell \in [1,c+d-k]$,
\begin{align*}
f(\ell)+f(\ell+1)&=C(\ell)+C(\ell+1)+(-1)^k(D(\ell)+D(\ell+1)\\
&=\sum_{r=0}^{k-1}(-1)^r\binom{c-1-r}{k-1-r}\binom{c+d-k}{r}\binom{p+k-r-1}{c-\ell-r}\\
&\quad+(-1)^k \sum_{j=0}^{k-1}(-1)^j\binom{d-1-j}{k-1-j}\binom{c+d-k}{j}\binom{p+k-1-j}{\ell+k-1-c-j}\\
&=\sum_{r=0}^{k-1}(-1)^r\binom{c-1-r}{k-1-r}\binom{c+d-k}{r}\binom{p+k-r-1}{c-\ell-r}\\
&\quad+(-1)^k \sum_{j=0}^{k-1}(-1)^j\binom{d-1-j}{k-1-j}\binom{c+d-k}{j}\binom{p+k-1-j}{p+c-\ell}.
\end{align*}
So $\mod p$
\begin{align*}
f(\ell)+f(\ell+1)&\equiv\sum_{r=0}^{k-1}(-1)^r\binom{c-1-r}{k-1-r}\binom{c+d-k}{r}\binom{k-1-r}{k+\ell-c-1}\\
&\qquad+(-1)^k\sum_{r=0}^{k-1}(-1)^r\binom{d-1-r}{k-1-r}\binom{c+d-k}{r}\binom{k-1-r}{c-\ell}\\
&=F_{c,d,k}(\ell)+(-1)^k G_{c,d,k}(\ell)\\
&=0
\end{align*}
by Lemma~\ref{FandG}.  Thus $f(\ell) \equiv (-1)^{\ell-1} f(1) \mod p$ for every integer $\ell \in [1,c+d-k]$.

\end{proof}

\section{Other Possible Results}\label{Others}
Theorem~\ref{Main} seems to be the first in a family of such results.  For example, a slightly more complicated result seems to hold when identifying a generator for the component $V_{\lambda_k}$ of $V_{p+c} \otimes V_{p+d}$ or for the component of $V_{2p-\lambda_k}$ in $V_{2 p+c}\otimes V_{2p+d}$ where $p \geq 5$.


\begin{thebibliography}{9}
\bibitem{Alperin}
J. L. Alperin,
Local representation theory, Cambridge Studies in Advanced Mathematics \textbf{11}, Cambridge University Press, Cambridge, 1986.
\bibitem{B2011_0}
M. J. J. Barry,
\emph{Decomposing Tensor Products and Exterior and Symmetric Squares},  J. Group Theory \textbf{14} (2011), 59--82.
\bibitem{B2017}
M. J. J. Barry,
\emph{Generators for decompositions of tensor products of modules associated with standard Jordan partitions}, Comm. Algebra \textbf{45} (2017), no. 4, 1819--1824.
\bibitem{GKP}
R. L. Graham, D. E. Knuth, and O. Patashnik,
\emph{Concrete Mathematics}, Second Edition,
Addison-Wesley, Reading MA, 1994.
\bibitem{Norman2008}
C. W. Norman,
\emph{On Jordan bases for the tensor product and Kronecker sum and their elementary divisors over fields of prime characteristic},
Linear and Multilinear Algebra \textbf{56} (2008), No. 4, 415--451.
\bibitem{Renaud1979}
J.-C. Renaud,
\emph{The decomposition of products in the modular representation ring of a cyclic group of prime power order}, J. Algebra \textbf{58} (1979), 1--11.
\end{thebibliography}
\end{document}